\documentclass[a4paper]{article}
\usepackage[english]{babel}
\usepackage{inputenc}
\usepackage{amsmath,latexsym,amsthm}
\usepackage{graphicx}
\usepackage{amssymb}
\begin{document}
\begin{flushright}
{\sf submitted to Uzbek  Mathematical Journal}

\end{flushright}
\sloppy

\begin{center}
\textbf{\large  Direct and inverse source problems for two-term time-fractional diffusion equation with Hilfer derivative}\\
\textbf{ Salakhitdinov M.S., Karimov E.T.}
\end{center}

\makeatletter
\renewcommand{\@evenhead}{\vbox{\thepage \hfil {\it  Salakhitdinov M.S., Karimov E.T.}   \hrule }}
\renewcommand{\@oddhead}{\vbox{\hfill
{\it Direct and inverse source problems for two-term time-fractional ...}\hfill
\thepage\\ \hrule}} \makeatother

\begin{abstract}

In this paper, we investigate direct and inverse source problems for the diffusion equation with two-term generalized fractional derivative (Hilfer derivative) in a rectangular domain. Using spectral expansion method, we derive two-term fractional differential equation together with appropriate initial condition (Cauchy problem). Based on solution of that Cauchy prob\-lem, we represent solution of formulated problems as a combination of sinus and multinomial Mittag-Leffler function of two variables. Imposing certain conditions to the given data, we prove uniform convergence of certain infinite series.
\end{abstract}

\textbf{MSC 2010:} 35R11, 26A33

\smallskip

\textbf{Keywords:} Hilfer derivative, multi-term fractional diffusion equation, inverse source problem, multivariate Mittag-Leffler function.

\section{Inverse source problem.}

Consider the following diffusion equation with Hilfer derivatives
\[
D_{0t}^{\alpha_1,\beta_1}u(t,x)+\mu D_{0t}^{\alpha_2,\beta_2}u(t,x)-u_{xx}(t,x)=g(x)\tag{1}
\]
in a rectangular domain $\Omega=\{(t,x):\,0<t<T,\,0<x<1\}$. Here $\alpha_i,\,\beta_i\,(i=1,2),\,\mu,\,T$ are given real numbers such that $T>0,\,0<\alpha_2<\alpha_1<1,\,\,0\leq\beta_i\leq 1$,
\[
D_{0t}^{\alpha,\beta}f(t)=\left(I_{0t}^{\beta(n-\alpha)}\frac{d^n}{dt^n}\left(I_{0t}^{(1-\beta)(n-\alpha)}f\right)\right)(t),\,\,t>0 \tag{2}
\]
is a Hilfer fractional derivative of the order $\alpha$ and type $\beta$ with $n-1<\alpha\leq n\in \mathbb{N}$ (see [1]),
\[
I_{0t}^\alpha f(t)=\frac{1}{\Gamma(\alpha)}\int\limits_0^t (t-z)^{\alpha-1}f(z)dz,\,\,t>0 \tag{3}
\]
is the Riemann-Liouville fractional integral of order $\alpha>0$ such that $I_{0t}^0 f(t)=f(t)$ (see, for example [2]).

We note that in case of $\beta=0$, Hilfer derivative (2) coincides with the Riemann-Liouville derivative and in case of $\beta=1$ with the Caputo derivative [1].

{\bf Problem 1} is to find a pair of functions $\{u(t,x),\,g(x)\}$ from the class of functions
\[
\begin{array}{l}

W_1=\{u(t,x):\,u\in C_{-1}[0,T],\,I_{0t}^{(1-\beta_1)(1-\alpha_1)}u\in C_{-1}^1[0,T],\\
u\in C[0,1],\,u_{xx}\in C(0,1);\,\,g(x)\in C[0,1]\},
\end{array}
\] 
which satisfies equation (1) in $\Omega$ together with boundary conditions
\[
u(t,0)=u(t,1)=0,\,\,0\leq t\leq T \tag{4}
\]
and initial conditions
\[
\lim\limits_{t\rightarrow +0} I_{0t}^{(1-\beta_i)(1-\alpha_i)}u(t,x)=0,\,\,0\leq x\leq 1\,\,(i=1,2), \tag{5}
\]
and also over-determining condition
\[
u(T,x)=\varphi(x),\,\,0\leq x \leq 1.\tag{6}
\]
Here $\varphi(x)$ is a given function such that $\varphi(0)=\varphi(1)=0$.

We remind the definition of the space $C_\alpha^m$.

{\bf Definition.} A function $f(x),\,\,x>0$, is said to be in the space $C_\alpha^m$, $m\in \mathbb{N}_0=\mathbb{N}\cup \{0\}$, $\alpha\in \mathbb{R}$ if there exists a real number $p,\,p>\alpha$, such that $f^{(m)}(x)=x^pf_1(x)$ with a function $f_1\in C[0,\infty)$, where $C_\alpha^0=C_\alpha$.

Hilfer derivative has appeared in the theoretical modeling of a broadband dielectric relaxation spectroscopy for glasses [3]. Some properties and appli\-cations of this derivative were studied in [4-7]. We will use in this work result of the work [8], where explicit solution of the modified Cauchy problem for multi-term fractional differential equation with Hilfer derivatives found by operational method.

Separation of variables lead us to the following spectral problem \[
X''(x)-\lambda X(x)=0,\,\,X(0)=X(1)=0,
\]
whose eigenvalues are $\lambda_k=(k\pi)^2,\,k\in \mathbb{Z}$ and corresponding eigenfunctions $\{X_k=\sin k\pi x\}$. Since this system forms complete orthogonal basis, we can expand solution of the problem 1 by the following series:
\[
u(t,x)=\sum\limits_{k=1}^\infty U_k(t)\sin k\pi x,\,0\leq t\leq T \tag{7}
\]  
\[
g(x)=\sum\limits_{k=1}^\infty g_k\sin k\pi x,\tag{8}
\]
where
\[
U_k(t)=2\int\limits_0^1 u(t,x)\sin k\pi x dx,\,\,k=1,2,... \tag{9}
\]
\[
g_k=2\int\limits_0^1 g(x)\sin k\pi x dx,\,\,k=1,2,... \tag{10}
\]

Substituting (7)-(8) into (1), we get
\[
D_{0t}^{\alpha_1,\beta_1}U_k(t)+\mu D_{0t}^{\alpha_2,\beta_2}U_k(t)+(k\pi)^2U_k(t)=g_k. \tag{11}
\]
Initial condition (5) gives us
\[
\lim\limits_{t\rightarrow +0} I_{0t}^{(1-\beta_i)(1-\alpha_i)}U_k(t)=0\,\,\,(i=1,2). \tag{12}
\]
According to [8], solution of (11) together with (12) has a form
\[
U_k(t)=g_k\int\limits_0^tz^{\alpha_1-1}E_{(\alpha_1-\alpha_2,\alpha_1),\alpha_1}\left(-\mu z^{\alpha_1-\alpha_2},-(k\pi)^2z^{\alpha_1}\right)dz, \tag{13}
\]
where
\[
E_{(\alpha-\beta,\alpha),\alpha}\left(x,y\right)=\sum\limits_{n=0}^\infty\sum\limits_{i=0}^{n}\frac{n!}{i!(n-i)!}\frac{x^i y^{n-i}}{\Gamma(\rho+\alpha n-\beta i)}\tag{14}
\]
is a particular case (in two variables) of multi-variate Mittag-Leffler function [9] with $\alpha,\beta,\rho>0$.

Applying the following formula (see, for example [10])
\[
\int\limits_0^t z^{\alpha-1}E_{(\alpha-\beta,\alpha),\alpha}\left(m_1 z^{\alpha-\beta},m_2z^{\alpha}\right)dz=t^{\alpha}E_{(\alpha-\beta,\alpha),\alpha+1}\left(m_1t^{\alpha-\beta},m_2t^{\alpha}\right)
\]
from (13) we obtain
\[
U_k(t)=g_k\cdot t^{\alpha_1-1}E_{(\alpha_1-\alpha_2,\alpha_1),\alpha_1+1}\left(-\mu t^{\alpha_1-\alpha_2},-(k\pi)^2t^{\alpha_1}\right).\tag{15}
\]
In order to define $g_k$ we use over-determined condition (6), which will take a form
\[
U_k(T)=\varphi_k,\tag{16}
\]
where
\[
\varphi_k=2\int\limits_0^1\varphi(x)\sin k\pi x dx,\,\,k=1,2,...
\]
From (15)-(16) we deduce
\[
g_k=\frac{\varphi_k}{T^{\alpha_1-1}E_{(\alpha_1-\alpha_2,\alpha_1),\alpha_1+1}\left(-\mu T^{\alpha_1-\alpha_2},-(k\pi)^2 T^{\alpha_1}\right)} 
\]
with
\[
{T^{\alpha_1-1}E_{(\alpha_1-\alpha_2,\alpha_1),\alpha_1+1}\left(-\mu T^{\alpha_1-\alpha_2},-(k\pi)^2 T^{\alpha_1}\right)}\neq 0. \tag{17}
\]
Finally, $U_k(t)$ will have a form
\[
U_k(t)=\frac{t^{\alpha_1-1}E_{(\alpha_1-\alpha_2,\alpha_1),\alpha_1+1}\left(-\mu t^{\alpha_1-\alpha_2},-(k\pi)^2 t^{\alpha_1}\right)}{T^{\alpha_1-1}E_{(\alpha_1-\alpha_2,\alpha_1),\alpha_1+1}\left(-\mu T^{\alpha_1-\alpha_2},-(k\pi)^2 T^{\alpha_1}\right)}\varphi_k. \tag{18}
\]
Now, let us estimate $U_k(t)$. Since
\[
\frac{t^{\alpha_1-1}E_{(\alpha_1-\alpha_2,\alpha_1),\alpha_1+1}\left(-\mu t^{\alpha_1-\alpha_2},-(k\pi)^2 t^{\alpha_1}\right)}{T^{\alpha_1-1}E_{(\alpha_1-\alpha_2,\alpha_1),\alpha_1+1}\left(-\mu T^{\alpha_1-\alpha_2},-(k\pi)^2 T^{\alpha_1}\right)}\leq C_1\,\,(C_1=const>0),
\]
we have
\[
\left|U_k(t)\right|\leq \frac{C_1}{(k\pi)^2}\left|\varphi_k^{(2)}\right|, \tag{19}
\]
where
\[
\varphi_k^{(2)}=-2\int\limits_0^1\varphi''(x)\sin k\pi x dx.
\]
We need more "strong estimate" \, for $U_k(t)$ in order to provide convergence of infinite series corresponding for $u_{xx}(t,x)$. Precisely, 
\[
\left|U_k(t)\right|\leq \frac{C_2}{(k\pi)^3}\left|\varphi_k^{(3)}\right|, \tag{20}
\]
where $C_2$ is a positive constant and
\[
\varphi_k^{(3)}=-2\int\limits_0^1 \varphi'''(x)\cos k\pi x dx.
\]
We have to impose more condition to the given function $\varphi(x)$ in order to guarantee uniform convergence of the following series
\[
u_{xx}(t,x)=\sum\limits_{k=1}^\infty U_k(t)(k\pi)^2 \sin k\pi x. \tag{21}
\]
In fact, considering (20), we have
\[
\left|u_{xx}(t,x)\right|\leq \sum\limits_{k=1}^\infty \left|U_k(t)\right| (k\pi)^2\leq \sum\limits_{k=1}^\infty \frac{C_2}{k\pi}\left|\varphi_k^{(3)}\right|.
\]
If we use $2ab\leq a^2+b^2$, we have
\[
\left|u_{xx}(t,x)\right|\leq \sum\limits_{k=1}^\infty \left( \frac{C_2^2}{4(k\pi)^2}+\left|\varphi_k^{(3)}\right|^2\right).
\]
Due to $\sum\limits_{k=1}^\infty \frac{1}{(k\pi)^2}=1/6$ and $\sum\limits_{k=1}^\infty |\varphi|^2\leq \left\|\varphi\right\|_{L_2(0,1)}$, we assume that
\[
\varphi(x)\in C^2[0,1],\,\varphi'''(x)\in L_2(0,1),\,\varphi(0)=\varphi(1)=\varphi''(0)=\varphi''(1)=0,\tag{22}
\]
then by Weierstrass M-test theorem, we can conclude that series (21) uni\-form\-ly converges.

Proof of the uniform convergence of series corresponding to the functions $u(t,x)$, $D_{0t}^{\alpha_i,\beta_i}u(t,x)\,(i=1,2)$ and $g(x)$ can be done by similar way consi\-dering (2), (3), (17), (19), but with less conditions to the given function $\varphi(x)$.

The uniqueness of the solution of the problem 1, can be obtained based on the completeness of the system $\{\sin k\pi x,\,k=1,2,...\}$ in $L_2$. In fact, if we consider corresponding homogeneous problem, i.e. $\varphi(x)=0$, from (18) we will get $U_k(t)\equiv 0$, which implies
\[
\int\limits_0^1 u(t,x)\sin k \pi x dx=0,\,\,0\leq t\leq T.
\]

Due to the completeness of the system $\{\sin k\pi x,\,k=1,2,...\}$ in $L_2$, we will get $u(t,x)\equiv 0$ in $\bar \Omega$.

We proved the following theorem:

{\bf Theorem 1.} If conditions (17) and (22) are valid, then the problem 1 is uniquely solvable and solution is represented by
\[
u(t,x)=\sum\limits_{k=1}^\infty  \frac{t^{\alpha_1-1}E_{(\alpha_1-\alpha_2,\alpha_1),\alpha_1+1}\left(-\mu t^{\alpha_1-\alpha_2},-(k\pi)^2 t^{\alpha_1}\right)}{T^{\alpha_1-1}E_{(\alpha_1-\alpha_2,\alpha_1),\alpha_1+1}\left(-\mu T^{\alpha_1-\alpha_2},-(k\pi)^2 T^{\alpha_1}\right)}\varphi_k\, sin k\pi x,
\] 
\[
g(x)=\sum\limits_{k=1}^\infty \frac{1}{T^{\alpha_1-1}E_{(\alpha_1-\alpha_2,\alpha_1),\alpha_1+1}\left(-\mu T^{\alpha_1-\alpha_2},-(k\pi)^2 T^{\alpha_1}\right)}\varphi_k\, sin k\pi x.
\]

\section{Direct problem.}

Now let us consider the following direct problem.

\textbf{Problem 2.} To find a solution of the equation
\[
D_{0t}^{\alpha_1,\beta_1}u(t,x)+\mu D_{0t}^{\alpha_2,\beta_2}u(t,x)-u_{xx}(t,x)=\bar g(t,x)\tag{23}
\]
from the class of functions
\[
\begin{array}{l}

W_2=\{u(t,x):\,u\in C_{-1}[0,T],\,I_{0t}^{(1-\beta_1)(1-\alpha_1)}u\in C_{-1}^1[0,T],\\
u\in C[0,1],\,u_{xx}\in C(0,1)\},
\end{array}
\]
satisfying conditions (4) and (5).

Here $\bar g(t,x)$ is a given function.

Similarly as in the case of problem 1, we search solution in the form of
\[
u(t,x)=\sum\limits_{k=1}^\infty \bar U_k(t)\sin k\pi x,\,0\leq t\leq T.\tag{24}
\]
Substituting (24) into (23) we get
\[
D_{0t}^{\alpha_1,\beta_1}\bar{U}_k(t)+\mu D_{0t}^{\alpha_2,\beta_2}\bar{U}_k(t)+(k\pi)^2\bar{U}_k(t)=\bar g_k(t), \tag{25}
\]
where
\[
\bar g_k(t)=2\int\limits_0^1 \bar g(t,x)\sin k \pi x dx.
\]

Solution of (25) satisfying initial condition
\[
\lim\limits_{t\rightarrow +0} I_{0t}^{(1-\beta_1)(1-\alpha_1)}\bar U_k(t)=0
\]
has a form
\[
\bar{U}_k(t)=\int\limits_0^tz^{\alpha_1-1}E_{(\alpha_1-\alpha_2,\alpha_1),\alpha_1}\left(-\mu z^{\alpha_1-\alpha_2},-(k\pi)^2z^{\alpha_1}\right)\bar g_k(t-z)dz. \tag{26}
\]
For the estimation of $\bar U_k(t)$ we use two different estimation of the function (14). Precisely, first of them is
\[
\left|E_{(\alpha-\beta,\alpha),\alpha}(x,y)\right|\leq \frac{C_3}{1+|x|}\tag{27}
\]
with $C_3=const>0$, which is proved in [9]. Another one is
\[
\left|E_{(\alpha-\beta,\alpha),\alpha}(x,y)\right|\leq\frac{C_4}{1+|x+y|},\tag{28}
\]
with $C_4=constant>0$, which has the following additional condition to the fractional orders (see [10], lemma 1.3)
\[
\Gamma(\rho+n(\alpha-\beta)+k\beta)>\Gamma(\rho+n(\alpha-\beta)),\,n,k\in\mathbb{N},\,n\geq k. \tag{29}
\]
If we use estimation (27), we get
\[
\left|\bar U_k(t)\right|\leq\frac{C_5}{(k\pi)^2}\left|\bar g_k^{(2)}(t)\right|,
\]
\[
\bar g_k^{(2)}(t)=-2\int\limits_0^1 \frac{\partial ^2\bar g(t,x)}{\partial x^2}\sin k \pi x dx.
\]
As we mentioned in previous case, we need another estimation for the $\bar U_k(t)$ in order to guarantee uniform convergence of infinite series corresponding to the function $u_{xx}(t,x)$, namely
\[
\left|\bar U_k(t)\right|\leq \frac{C_6}{(k\pi)^3}\left|\bar g_k^{(3)}(t)\right|,\tag{30}
\]
\[
\bar g_k^{(3)}(t)=-2\int\limits_0^1 \frac{\partial ^3\bar g(t,x)}{\partial x^3}\cos k \pi x dx.
\]
We impose the following conditions to the $\bar g(t,x)$:
\[
\begin{array}{l}
\frac{\partial ^2\bar g(t,x)}{\partial x^2}\in C[0,1],\,\frac{\partial ^3\bar g(t,x)}{\partial x^3}\in L_2(0,1),\\
\bar g(t,0)=\bar g(t,1)=0,\,\,\left.\frac{\partial ^2\bar g(t,x)}{\partial x^2}\right|_{x=0}=\left.\frac{\partial ^2\bar g(t,x)}{\partial x^2}\right|_{x=1}=0 \tag{31}
\end{array}
\]
in order to get
\[
\left|u_{xx}(t,x)\right|\leq C_7+\left\|\bar g_k^{(3)}(t)\right\|, \tag{32}
\]
where $C_7$ is a positive constant.

Now, if we use estimation (28), we obtain
\[
\left|\bar U_k(t)\right|\leq \frac{C_8}{(k\pi)^2}\left|\bar g_k^{(1)(t)}\right|,
\]
\[
\bar g_k^{(1)}(t)=2\int\limits_0^1\frac{\partial \bar g(t,x)}{\partial x}\cos k\pi x dx.
\]
Consequently, in order to provide the uniform convergence of series
\[
u_{xx}(t,x)=\sum\limits_{k=1}^\infty \bar U_k(t)(k\pi)^2\sin k\pi x \tag{33}
\]
we impose the following condition to the given function $\bar g(t,x)$:
\[
\frac{\partial \bar g(t,x)}{\partial x}\in C[0,1],\,\frac{\partial^2 \bar g(t,x)}{\partial x^2}\in L_2(0,1),\,\bar g(t,0)=\bar g(t,1)=0, \tag{34}
\]
which yields
\[
\left|u_{xx}(t,x)\right|\leq C_9+\left\|\bar g_k^{(2)}(t)\right\|.
\]
Using Weierstrass M-test theorem, one can easily prove the uniform conver\-gence of (33). The uniqueness of the solution for the problem 2 can be proved similarly to the proof of the problem 1.

Hence, we proved the following theorems:

\textbf{Theorem 2.} If condition (31) is valid, then problem 2 is uniquely solvable and solution is represented by
\[
\begin{array}{l}
u(t,x)=\\
\sum\limits_{k=1}^\infty \sin k\pi x \int\limits_0^tz^{\alpha_1-1}E_{(\alpha_1-\alpha_2,\alpha_1),\alpha_1}\left(-\mu z^{\alpha_1-\alpha_2},-(k\pi)^2z^{\alpha_1}\right)\bar g_k(t-z)dz. \tag{35}
\end{array}
\]

\textbf{Theorem 3.} If condition (29) and (34) are valid, then problem 2 is uniquely solvable and solution is represented by (35).

\medskip

\textbf{References}

{\small
\begin{enumerate}
\item R. Hilfer, Applications of Fractional Calculus in Physics. World Scientific, Singapore (2000).

\item I. Podlubny, Fractional Differential Equations, in: Mathematics in Science and Engineering, vol. 198, Acad. Press, 1999.

\item R. Hilfer, Experimental evidence for fractional time evolution in glass forming materials. Chem. Phys. 284 (2002), 399-408.

\item R. Hilfer, Y. Luchko, ?Z. Tomovski, Operational method for the solution of fractional differential equations with generalized Riemann-Liouville fractional derivatives. Fract. Calc. Appl. Anal. 12, No 3 (2009), 299-318.

\item  K.M. Furati, M.D. Kassim, N.e.-Tatar, Existence and uniqueness for a problem involving Hilfer fractional derivative. Comp. Math. Appl. 64, No 6 (2012), 1616–1626.

\item  K.M. Furati, O.S. Iyiola, M. Kirane, An inverse problem for a generalised fractional diffusion, Appl. Math. Comput. 249 (2014) 24–31.

\item S.A.Malik, S.Aziz. An inverse source problem for a two parameter anomalous diffusion equation with nonlocal boundary conditions. Computer and Mathe\-matics with Applications, 2017,\\ http://dx/doi.org/10.1016/j.camwa.2017.03.019

\item M.-Ha. Kim, G.Chol-Ri, Chol O.H. Operational method for solving multi-term fractional differential equations with  the generalized fractional deriva\-tives. Fractional Calculus and Applied Analysis. 2014. vol.17. No 1, pp.79-95.

\item  Luchko Y. and  Gorenflo R.  An Operational Method for Solving Fractional Differential Equations with the Caputo Derivatives, Acta Math.Vietnamica, N 2, 1999, 207-233.

\item E.T. Karimov, S. Kerbal, N. Al-Salti. Inverse source problem for multi-term fractional mixed type equation. M.Ruzhansky et al. (eds.), Advances in Real and Complex Analysis with Applications, Trends in Mathematics. 2017, pp. 289-301

\end{enumerate}
}

\begin{tabular}{p{7cm}l}
Institute of Mathematics &  \\
named after V.I.Romanovsky
\end{tabular}

\end{document}